\documentclass[amssymb,10pt]{article}

\usepackage{amssymb}
\usepackage{amsthm}
\usepackage{indentfirst}
\usepackage{amsmath}

cm \hoffset = -2.5 true cm

\newtheorem{theoremalph}{Theorem}

\newtheorem*{Theorem A}{Theorem A}

\newtheorem*{Conjecture}{Conjecture}

\newtheorem*{Problem}{Problem}

\def\DD{{\mathbb D}}

 \def\NN{{\mathbb N}} 

 \def\RR{{\mathbb R}}

 \def\ZZ{{\mathbb Z}}

\def\Si{\Sigma}
\def\La{\Lambda}

\def\De{\Delta}

\def\Diff{\hbox{Diff} }

\def\cB{{\cal B}}    
\def\cC{{\cal C}}    \def\cU{{\cal U}}
\def\cD{{\cal D}}   \def\cP{{\cal P}} \def\cV{{\cal V}}
    
\def\cF{{\cal F}}  \def\cL{{\cal L}} \def\cR{{\cal R}} \def\cX{{\cal X}}

\def\Int{\operatorname{Int}}
\def\dim{\operatorname{dim}}

\newtheorem{theo}{Theorem}

\newtheorem{clai}{Claim}

\newtheorem{lemm}{Lemma}[section]
\newtheorem{coro}[lemm]{Corollary}
\newtheorem{defi}[lemm]{Definition}
\newtheorem{prop}[lemm]{Proposition}
\newtheorem{conj}{Conjecture}

\newtheorem{prob}{Problem}

\newtheorem{ques}{Question}

\newtheorem{rema}[lemm]{Remark}

\newenvironment{demo}[1][Proof]{{\bf #1~: }}{\hfill$\Box$\medskip}

\begin{document}

\title{On the existence of attractors}
\author{
Christian Bonatti \and Ming Li \and  Dawei Yang\footnote{This work
has been done during the stays of Li Ming and Yang Dawei at the
IMB, Universit\'e de Bourgogne and we thank the IMB for its warm
hospitality. M. Li is supported by a post doctoral grant of the
R\'egion Bourgogne, and D. Yang is supported by CSC of Chinese
Education Ministry.} }

\maketitle

\begin{abstract} On every compact $3$-manifold, we build a non-empty open set $\cU$ of $\Diff^1(M)$ such that, for every $r\geq 1$, every $C^r$-generic diffeomorphism $f\in\cU\cap \Diff^r(M)$ has no topological attractors. On higher dimensional manifolds, one may require that $f$ has neither topological attractors nor topological repellers. Our examples have finitely many \emph{quasi attractors}. For flows, we may require that these quasi attractors contain singular points.  Finally we discuss alternative definitions of attractors which may be better adapted to generic dynamics.

\end{abstract}

\section{Introduction}
 The aim of dynamical systems is to describe  the asymptotic behavior of the orbits when the time tends to infinity.
For simple dynamical systems, the behavior of the orbits looks
like the gradient flow  of a Morse function: most of the orbits
tend to a sink, and the union of the basins of the sink is a dense
open set in the ambient manifold.

However, many dynamical systems present a more complicated
behavior and many orbits do not tend to periodic orbits;  their
$\omega$-limit set may be chaotic. In the sixties and seventies,
many people tried to give a definition of attracting sets,
allowing to describe most of the possible behaviors of dynamical
systems. An attractor $\Lambda$ of a diffeomorphism $f$ needs to
satisfy two kinds of properties:
\begin{itemize}
\item it attracts ``many orbits''. According to the authors, this
means: the basin of $\Lambda$ contains a neighborhood of
$\Lambda$, an open set,  a residual subset of an open set, a set
with positive Lebesgue measure, \dots.

\item it is indecomposable, that is, it cannot split into the
union of smaller attractors; many notions of indecomposability are
used:
 transitivity (generic orbits of the attractor are dense in the attractor), chain recurrence (for every $\delta>0$, one can go from any point of the attractor to any point of the attractor by $\delta$-pseudo orbits inside the attractor), uniqueness of the SRB measure,\dots

\end{itemize}
None of these notions can cover all the possible behaviors of
dynamical systems. For every notion of (indecomposable)
attractors, one can find examples of dynamical systems without
attractors. \footnote{If one removes the indecomposability
hypothesis,  Conley shows that attractors  exist for any
homeomorphism of compact metric space. More precisely, given any
points $x, y$, either one can join $x$ to  $y$ by
\emph{$\delta$-pseudo orbits}, for every $\delta>0$, or there is
an \emph{ attracting region $U$} containing $x$ but not $y$. Hence
the dynamics admits attracting regions, or the chain recurrent set
is the whole space. Conley calls \emph{attractors} the maximal
invariant sets in the attracting regions.
 The attractors in Conley theory are not  assumed to be indecomposable: an attractor can contain smaller attractors.}

A natural idea for bypassing this difficulty is to restrict the
study to generic dynamical systems, in order to avoid the most
pathological and fragile behaviors.  A property is $C^r$-generic
if it holds on a residual subset of the space of $C^r$
diffeomorphisms $\Diff^r(M)$ endowed with the $C^r$ topology.

This viewpoint has been considered very early by Smale and Thom,
with the hope that generic dynamical systems would have a simple
behavior. For instance one can read in \cite[Chapter 4.1 B]{T}:
\emph{Il n'est pas certain qu'un champ $X$ donn\'e dans $M$
pr\'esente des attracteurs, a fortiori des attracteurs
structurellement stables. Toutefois, selon certaines id\'ees
r\'ecentes de Smale, si la vari\'et\'e $M$ est compacte, presque
tout champ pr\'esenterait un nombre fini d'attracteurs isol\'ement
structurellement stables;(\dots)}\footnote{\emph{``It is not clear
if a given vector field in $M$ has an attractor, \emph{a fortiori}
a structurally stable attractor; however, according to recent
ideas by Smale, if the manifold is compact, almost all vector
fields would admit finitely many attractors, each of them
structurally stable;(\dots)."}}. Thom's idea was renewed and
formalized in 1975 as \emph{Thom's conjecture} by Palis and Pugh
{\cite[Problem 26]{PP}}:\emph{ There is a dense open set in ${\rm
Diff}^r(M)$ such that for almost every point $x\in M$, the
$\omega$-limit set $\omega(x)$ is  a topological attractor, and
each attractor is topologically stable. }

After thirty years of progress in the field, this conjecture can
look naively optimistic. Indeed, Thom's original idea was
disproved in most of its aspects: finiteness, stability.
\begin{itemize}
\item there are open sets of systems without structurally stable
attractors, as the robustly transitive non-hyperbolic
diffeomorphisms built by Shub in \cite{Sh}; \item there are
$C^r$-locally generic diffeomorphisms having infinitely many sinks
(see \cite{N1,N2} for $r\geq 2$ and \cite{BD2} for $r=1$).
\end{itemize}

However, the existence of at least one attractor remains an open
question. In this paper, we will give a negative answer to this
question, showing that the usual notion of topological attractor
is too strong and not adapted to generic dynamical systems. Let us
be now somewhat more precise.

\vskip 5mm A \emph{topological attractor} of a diffeomorphism
$f\colon M\to M$ is a compact subset $\Lambda\subset M$ with the
following properties
\begin{itemize}
\item $\La$ is \emph{invariant} (i.e. $f(\Lambda)=\Lambda$);

\item $\La$ admits a compact neighborhood $U$ which is \emph{an
attracting region} (i.e. the image $f(U)$ is contained in the
interior of $U$) such that all the orbits in $U$ converge to
$\La$: $\La=\bigcap_{n\in\NN} f^n(U)$;

\item $\La$ is  \emph{transitive} (i.e. the positive orbits of
generic points in $\Lambda$ are dense in $\Lambda$), or at least
\emph{chain recurrent}\footnote{some authors say ``\emph{chain
transitive}''.} (i.e. any two points in $\La$ can be joined by
$\varepsilon$-pseudo orbits, for all $\varepsilon>0$). We will
speak on \emph{transitive topological attractor}, and on
\emph{chain recurrent topological attractor}, if we need to
emphasize the kind of indecomposability we require.
\end{itemize}
A \emph{topological repeller} of $f$ is, by definition,  a
topological attractor of $f^{-1}$.
 In 2004 \cite{BC} and in 2005 \cite[Problem 10.35]{BDV} still asked:
\begin{ques} Do $C^1$-generic diffeomorphisms
admit at least one (chain recurrent or transitive \footnote{Chain
recurrent topological attractors of $C^1$-generic diffeomorphisms
are homoclinic classes, hence are transitive.}) topological
attractor? Does the union of the basins of the topological
attractors cover a dense open subset of the manifold?
\end{ques}

The answer is ``no":  The fact that $C^0$-generic homeomorphisms
have no attractors is known from Hurley's work \cite{H}:  every
attracting region of a $C^0$-generic homeomorphism contains
infinitely many repelling regions and infinitely many disjoint
attracting regions. Theorems~\ref{Thm:3Dnoattractor}
and~\ref{Thm:4Dnoattractorrepellers} show that, for every $r\geq
1$, the property of having  (at least) one topological attractor
is not $C^r$-generic.

Our results use the notion of \emph{quasi attractors}, introduced
by Hurley: a chain recurrence class of a homeomorphism is a
\emph{quasi attractor}\footnote{Some authors use the terminology
``weak attractor'' instead of quasi attractor.} if it is the
intersection of a sequence of attracting regions. A  \emph{quasi
repeller} for $h$ is a quasi attractor for $h^{-1}$.

\begin{theoremalph}\label{Thm:3Dnoattractor}
For every three-dimensional manifold $M^3$, there is a  non-empty
$C^1$-open subset ${\cal U}\subset {\rm Diff}^1(M^3)$ such that:
\begin{itemize}
\item there are hyperbolic periodic saddles $p_{1,f},\dots,
p_{k,f}$ varying continuously with $f\in\cal U$, whose chain
recurrence classes $\La_{1,f}, \dots,\La_{k,f}$ are the unique
quasi attractors of $f$; \item the set $\{f\in\cU, f \mbox{ has no
attractors}\}$ is $C^r$-residual in $\cU\cap {\rm Diff}^r(M)$ for
every $r\geq 1$.
\end{itemize}
\end{theoremalph}

The $C^r$-generic diffeomorphisms $f$ in the open set $\cU$ have
no attractors but infinitely many repellers. This motivates the
following problem:
\begin{Problem}
For three-dimensional manifold $M^3$, is there a dense (open and
dense, residual) set
 $\cD\subset {\rm Diff}^r(M)$, such that  any
$f\in\cD $  has neither attractors nor repellers?
\end{Problem}

\begin{theoremalph}\label{Thm:4Dnoattractorrepellers}
For every compact  manifold $M$, with $\dim M\geq 4$ there is a
non-empty open set ${\cal U}\subset {\rm Diff}^1(M)$ such that:
\begin{itemize}
\item there are hyperbolic periodic saddles $p_{1,f},\dots,
p_{k,f}$ varying continuously with $f\in\cal U$, whose chain
recurrence classes $\La_{1,f}, \dots,\La_{k,f}$ are the unique
quasi attractors of $f$; \item there are hyperbolic periodic
saddles $q_{1,f},\dots, q_{\ell,f}$ varying continuously with
$f\in\cal U$, whose chain recurrence classes $\Si_{1,f},
\dots,\Si_{\ell,f}$ are the unique quasi repellers of $f$; \item
the set $\{f\in\cU, f \mbox{ has neither attractors  nor
repellers}\}$ is $C^r$-residual in $\cU\cap {\rm Diff}^r(M)$ for
every $r\geq 1$.
\end{itemize}
\end{theoremalph}

Our results can be easily adapted for vector fields, building
locally generic vector fields having finitely many (non-singular)
quasi attractors but no attractors. However, one of the main
differences between diffeomorphisms and flows is the existence of
singularities, in particular when these singularities are not
isolated from the regular part of the limit set of the flow.

This phenomenon has first been suspected experimentally by Lorenz
\cite{Lo}, and then proved rigourously in \cite{Gu,ABS,GuW}, where
the authors exhibited, in dimension $3$,  a $C^1$-open set of
vector fields having a robust attractor containing infinitely many
periodic orbits accumulating on a saddle singularity. Their
construction (known as \emph{geometric model of Lorenz attractor})
leads to the notion of singular attractors, which have been
studied in extends on $3$-manifolds: for instance, if the presence
of a singularity inside the attractor prevents the usual
definition of hyperbolicity, robust singular attractors in
dimension $3$ always present a kind  of weak hyperbolicity called
{\it singular hyperbolicity}, see \cite{MPP1,MPP2}. In particular,
they satisfy the \emph{star condition}: $C^1$-robustly all the
periodic orbits are hyperbolic. \cite{LGW,ZGW,MM} show that in any
dimension, robust singular attractors satisfying the star
condition are singular hyperbolic. Recent examples \cite{BKR}
\cite{BLY} show that robust singular attractors may  satisfy
neither the star condition nor the singular hyperbolicity. However
even these new examples admit a strong stable direction, invariant
by the flow and dominated by a center-unstable bundle.

Hence it is natural to ask:
\begin{ques} Does every $C^1$-robust singular attractor admit a strong stable bundle?\end{ques}

Indeed, this question has been our first motivation for this work.
Before presenting our results, let us make a comment on this
question. Examples of (non-singular) robustly transitive attractor
whose flow does not admit any dominated splitting are already
known (just consider the suspension flows of robustly transitive
diffeomorphisms without invariant hyperbolic bundles in
\cite{BV}). For this reason one considers the linear Poincar\'e
flow on the normal bundle and this flow admits a dominated
splitting. However, the linear Poincar\'e flow is not defined on
the singularity: for this reason, it is not clear what kind of
hyperbolicity will satisfy the singular attractors.

Now we state our result for flows. Our construction can be adapted
in order to build a \emph{robust singular quasi attractor} whose
tangent bundle doesn't have any dominated splitting with respect
to the tangent flow.

\begin{theoremalph}\label{Thm:4Dvectorfield}
There is a non-empty open set ${\cal U}$ of the space
${\cX}^r(B^4)$ of $C^r$ vector fields on the $4$-ball,  such that:
\begin{itemize}
\item any $X\in\cal U$  is transverse to the boundary and entering
inside the ball; \item  any $X\in\cal U$ has a unique zero $0_X$
in $B^4$; one denotes by $\La_X$ the chain recurrence class of
$0_X$; \item any $X\in\cal U$ has a unique quasi attractor in
$B^4$ which is $\La_X$; \item the subset $\{X\in\cU,\ \La_X \mbox{
is not an attractor}\}$ is $C^r$-residual in $\cU$; \item for
$X\in \cU$, there is no dominated splitting for the tangent flow
of $X$ on $\La_X$.
\end{itemize}
\end{theoremalph}

\subsection{Organization of the paper}
Our main result is the construction in Section~\ref{ss.torus} of
an example of locally generic diffeomorphisms of the solid torus
$S^1\times \DD^2$, without attractors.

Putting the solid torus in a ball $B^3$, we get a model of an
attracting ball without attractors, which allows us, in
Section~\ref{ss.dim3} to replace the sinks of a gradient like
diffeomorphism by these attracting balls without attractors,
proving Theorem~\ref{Thm:3Dnoattractor}.

Multiplying this ball $B^3$ by a normal contraction, one gets in
Section~\ref{ss.dim4} an attracting ball $B^n$, for $n>3$, without
attractors and repellers. This section ends the proof of
Theorem~\ref{Thm:4Dnoattractorrepellers}.

Section~\ref{ss.vector} considers the case of vector fields and
shows that our construction in Section~\ref{ss.torus} leads to
locally generic vector fields  $X$ on $4$-manifolds having a
unique quasi attractor $\La_X$ and no attractors; furthermore
$\La_X$ is the chain recurrence class of a singularity of $X$.

Section~\ref{s.conclusion} concludes this paper by discussing
alternative notions of attractors which could be better adapted to
generic dynamical systems.

\section{Notations, definitions  and preliminaries}

\subsection{Disks and balls}
For every  $d\in\mathbb N$ and $r\in\RR$, we denote by $\DD^d(r)$
the  closed ball in ${\mathbb R}^d$ centered at $0$ and with
radius $r$, i.e., $\DD^d(r)=\{x\in{\mathbb R}^d:\|x\|\le r\}$. For
simplicity, we denote $\DD^d=\DD^d(1)$. Given a compact Riemannian
manifold $M$,  a point  $x\in M$, and a real number $\delta>0$, we
denote $B_\delta(x)=\{y\in M: ~d(x,y)\le\delta\}$, the compact
ball centered at $x$ and with radius $r$.

Recall that every orientation preserving diffeomorphism of $\DD^2$
is smoothly isotopic to the identity.

An \emph{essential disk in $S^1\times \DD^2$} is a embedding
$D\colon \DD^2\hookrightarrow S^1\times \DD^2$ whose boundary
$\partial D= D(\partial\DD^2)$ is contained in $\partial
\left(S^1\times \DD^2\right) =S^1\times S^1$, and is not homotopic
to a point in $S^1\times S^1$.

\subsection{Hyperbolicity, partial hyperbolicity, dominated
splitting} Let $f$ be a diffeomorphism of a manifold $M$ of
dimension $d$, $x$  a periodic point of $f$, and $\pi$  its
period. Let $\lambda_1\leq\lambda_2\leq\dots\leq \lambda_d$ be the
moduli of the eigenvalues of the differential $Df^\pi(x)$. The
point $x$ is \emph{hyperbolic} if $\lambda_i\neq 1$ for all
$i\in\{1,\dots,d\}$.  The point $x$ is \emph{sectionally area
expanding} (or \emph{sectionally expanding}) if
$$\lambda_i\lambda_j>1, \mbox{ for all }i,j\in\{1,\dots,d\}, i\neq
j.$$

We say a compact invariant set $\Lambda$ of $f$ is {\it
hyperbolic} if there are a $Df$-invariant splitting
$$TM|_\Lambda
=E^s\oplus E^u,$$ and constants $C>0$, $\lambda\in(0,1)$ such that
for any $x\in\Lambda$ and $n\in\mathbb N$
$$\|Df^{n}|_{E^s(x)}\|\le C\lambda^n,~~~\|Df^{-n}|_{E^u(x)}\|\le C\lambda^n.$$

The bundles $E^s$ and $E^u$ are called the \emph{stable} and
\emph{unstable bundle} of $\La$, respectively. They are always
continuous bundles, so that the dimensions $\dim{E^s(x)}$ and
$\dim{E^u(x)}$ are locally constant. If $\dim{E^s(x)}$ is
independent on $x\in\Lambda$, then we call $\dim{E^s}$  the index
of the hyperbolic set $\La$.

We say $\Lambda$ is a {\it basic set} if $\Lambda$ is an isolated
hyperbolic transitive set: there is an open neighborhood $U$ of
$\La$ such that $\La=\bigcap_{i\in\ZZ} f^i(U)$.

Given a compact invariant set $\Lambda$, a $Df$-invariant
splitting $TM|_{\La}= E_1\oplus E_2\oplus\dots \oplus E_k$ is
\emph{dominated}, and we denote $E_1\oplus_{_<}
E_2\oplus_{_<}\dots \oplus_{_<} E_k$, if the dimensions
$\dim(E_i)$ are constant over $\La$,  and if there are constants
$C>0$, $\lambda\in(0,1)$ such that for any $x\in\Lambda$,
$n\in\mathbb N$ and $i\in\{1,\dots,k-1\}$,  we have
$$\|Df^{n}|_{E_i(x)}\|\|Df^{-n}|_{E_{i+1}(f^n(x))}\|\le C\lambda^n.$$


A $Df$-invariant bundle $E$ is called \emph{(uniformly)
contracting} if there are constants $C>0$, $\lambda\in(0,1)$ such
that for any $x\in\Lambda$ and $n\in\mathbb N$, we have
$\|Df^{n}|_{E(x)}\|\le C\lambda^n$; it is called \emph{(uniformly)
expanding} if it is contracting for $f^{-1}$.

A dominated splitting $E_1\oplus_{_<} E_2\oplus_{_<}\dots
\oplus_{_<} E_k$ is \emph{partially hyperbolic} if $E_1$ is
uniformly contracting or $E_k$ is uniformly expanding.

\subsection{Cone fields associated to a dominated splitting}

Given a compact set $V\subset M$, a continuous (not necessary
invariant) bundle $F\subset T_V M$, and a positive number
$\alpha>0$, \emph{the cone field on $V$ associated to $F$ of size
$\alpha>0$} is
$$C^F_\alpha(x)=\{v\in T_x M:~\exists v_F\in F,v_{F^\perp}\in F^\perp,{~\rm s.t.~}v=v_F+v_{F^\perp}, |v_{F^\perp}|\le\alpha|v_F|\}$$
for $x\in V$, where $F^\perp$ is the orthogonal subbundle of $F$.

We say a \emph{cone field  $C^F_\alpha$ is strictly
$Df$-invariant}, if there is  $\beta\in(0,\alpha)$ such that, for
any $x\in V$ such that $f(x)\in V$, we have
$$Df(C^F_\alpha(x))\subset
C^F_\beta(f(x)).$$

If an invariant compact set $\Lambda$ has a dominated splitting
$T_{\Lambda}M=E\oplus_{_<} F$, then there is $\alpha_0>0$, such
that for any $\alpha\in(0,\alpha_0)$, there is $N\in\mathbb N$
such that the
 cone field $C^E_\alpha$
is strictly $Df^{-N}$-invariant and the cone field  $C^F_\alpha$
is strictly $Df^N$-invariant.

\subsection{Conley theory and quasi attractors}
Let $(X,d)$ be a compact metric space and $f\colon X\to X$ a
homeomorphism.

For any
 $x,y\in X$, we denote $x\dashv y$ if  for every $\epsilon>0$ there
is a $\epsilon$-pseudo orbit joining $x$ to $y$, that is: there
are $n>0$ and a sequence of points $\{x=x_0,x_1,\cdots, x_n=y\}$
verifying $d(f(x_i),x_{i+1})<\epsilon$ for $0\le i\le n-1$.

 We say that $x$ is chain recurrent if $x\dashv
x$, and we denote by $\cR(f)$ the set of chain recurrent points of
$f$, called the \emph{chain recurrent set} of $f$.

An invariant compact set $K$ of $X$ is \emph{chain recurrent} (or
\emph{chain transitive}) if every point $x\in K$ is chain
recurrent for the restriction $f|_K$: in other words,
$K=\cR(f|_K)$.

We say $x$ and $y$ are {\it chain equivalent} if $x\dashv y$ and
$y\dashv x$. The  chain equivalence is an equivalence relation on
$\cR(f)$. For any $x\in {\cR}(f)$, the equivalence class of $x$ is
called   {\it the chain recurrence class of $x$}, and denoted by
$C(x)$.


A \emph{quasi attractor} $\La$ is a chain transitive set which
admits a base of  neighborhood which are attracting regions (this
implies that $\La$ is a chain recurrence class).

\subsection{Plykin attractor}\label{ss.plykyn}

Let $\Diff_0^1(\DD^2,\Int(\DD^2))$ denote the space of orientation
preserving  $C^1$-embeddings  $\phi\colon \DD^2\to \Int(\DD^2)$.
Notice that the elements of $\Diff_0^1(\DD^2,\Int(\DD^2))$ are all
isotopic, in particular are isotopic to any linear contraction of
$\DD^2$.

In \cite{Pl} Plykin built a non-empty open subset $\cP\subset
\Diff_0^1(\DD^2,\Int(\DD^2))$  such that for any $\phi\in\cP$ the
chain recurrent set of $\phi$ consists in the union of a
non-trivial hyperbolic attractor $A_\phi$ and a finite set of
periodic sources.

We denote by $\cP_0\subset \cP$ the non-empty open subset of
diffeomorphisms such that the hyperbolic attractor $A_\phi$
contains a fixed  point $x_\phi$ which is an area expanding saddle
point:
 $$Det(D\phi (x_\phi))>1.$$

\subsection{Solenoid maps associated to a braid in $S^1\times \DD^2$}

A \emph{connected braid $\gamma $ of $S^1\times \DD^2$}  is (the
isotopy class of) an embedding of the circle $S^1$ in $S^1\times
\DD^2$, transverse to the fibers $\{\theta\}\times \DD^2$, for
$\theta\in S^1$. The projection $S^1\times \DD^2\to S^1$ induces
on $\gamma$ a finite covering of the circle; we denote by
$n_\gamma\neq 0$ the order of this finite cover.

For any  braid $\gamma$, we denote by $\cU_{\gamma}$  the
(non-empty) open subset of diffeomorphisms $f\colon S^1\times
\DD^2\hookrightarrow\Int(S^1\times \DD^2)$ such that
$f(S^1\times\{0\})$ is isotopic to the braid $\gamma$.

We call \emph{canonical solenoid maps} associated to a braid
$\gamma$ the maps built as follows:  denote $n=n_\gamma$;  we
choose a representative $\gamma\colon S^1\to S^1\times \DD^2$ of
the braid having the following form: $\gamma(t)=(n.t, z(t))$. We
fix $\delta>0$ such that
$$\begin{array}{rl}
\mbox{for all } t\in S^1,& d(z(t),\{n.t\}\times \partial\DD^2)>2\delta;\\
\mbox{ for any }t_1,t_2\in S^1,&\left( t_1\neq t_2\mbox{ and }
n.t_1=n.t_2\in S^1\right) \Rightarrow d(z(t_1), z(t_2))>2\delta.
\end{array}$$

Now the map $f_{\gamma,\delta}$ defined on $S^1\times \DD^2$ by
$f_{\gamma,\delta}(t,z)= (n.t, \delta.z + z(t))$ belongs to
$\cU_\gamma$ and is called a \emph{canonical solenoid map}
associated to a braid $\gamma$.

\subsection{Partially hyperbolic solenoid maps}

For every $\alpha>0$ we denote by $\cC_\alpha$ the cone field on
$S^1\times \DD^2$ defined by
$$\cC_\alpha(x)=\{u=(u_1,u_2)\in T_x(S^1\times \DD^2)=\RR\times
\RR^2 \mbox{ such that } |u_2|\leq \alpha |u_1|\}.$$

We denote by $\cU_\gamma^{part.hyp}$ the set of diffeomorphisms
$f\in\cU_\gamma$ such that there are $\alpha>0$ and $\ell\in
\NN\setminus\{0\}$  such that:
\begin{itemize}
\item  the cone field $\cC_\alpha$ is strictly invariant  under
$Df^\ell$; \item there is $\lambda >1$ such that for every $x\in
S^1\times \DD^2$ and every vector $u=(u_1,u_2)\in\cC_\alpha(x)$
one has:
$$|v_1|\geq \lambda |u_1| \mbox{, where } Df^\ell(u)=(v_1,v_2)\in
T_{f^\ell(x)}(S^1\times \DD^2).$$
\end{itemize}

The set $\cU_\gamma^{part.hyp}$ is a $C^1$-open subset of
$\cU_\gamma$. Moreover, one easily verifies:

\begin{lemm}\label{l.parthyp}
Let $f\in\cU_\gamma$ be of the form $(t,z)\mapsto (\gamma(t),
\varphi_t(z))$. Assume that:
\begin{itemize}
\item for every $t\in S^1$ one has
$$\left| \frac d{dt}\gamma(t)\right|>1;$$
\item for every $(t,z)\in S^1\times \DD^2$ one has
$$\left\| D_z(\varphi_t)\right\|<\left| \frac
d{dt}\gamma(t)\right|.$$
\end{itemize}
Then $f$ is partially hyperbolic; more precisely,
$f\in\cU_\gamma^{part.hyp}$.
\end{lemm}

As a direct consequence one gets:
\begin{coro}\label{c.canonical} For every braid $\gamma$ with $|n_\gamma|\geq 2$ every
canonical solenoid map  $f$ associated to $\gamma$ belongs to
$\cU_\gamma^{part.hyp}$.

In particular the open set $\cU_\gamma^{part.hyp}$ is non-empty.
\end{coro}

\begin{coro}\label{c.parthyp} Let $f\in \cU_\gamma$ satisfying the hypotheses of
Lemma~\ref{l.parthyp}, and $h\colon S^1\times \DD^2\to S^1\times
\DD^2$ be a diffeomorphism of the form $(t,z)\mapsto (t,h_t(z))$
where $h_t$ is an orientation preserving diffeomorphism of
$\DD^2$.

Then the map $g= h^{-1}f h$ belongs to $\cU_\gamma^{part.hyp}$.
\end{coro}


\subsection{Realization of a map $\varphi\in \Diff_0^1(\DD^2,\Int(\DD^2))$ by a solenoid map $f\in\cU_\gamma^{part.hyp}$ }

The aim of this section is to prove:
\begin{prop}\label{p.realization} Given any $\varphi\in {\rm Diff}_0^1(\DD^2,\Int(\DD^2))$ and
any braid $\gamma$ with $|n_\gamma|\geq 2$, there is a
diffeomorphism $f\in \cU_\gamma^{part.hyp}$ such that the disk
$\{0\}\times \DD^2$ is positively invariant (and normally
hyperbolic), and the restriction of $f$ to $\{0\}\times \DD^2$ is
$\varphi$.
\end{prop}
\begin{demo} We denote $n=n_\gamma$.  We
choose a representative $\gamma\colon S^1\to S^1\times \DD^2$,
$\gamma(t)=(n.t, z(t))$. Consider a canonical solenoid map
$f_{\gamma,\delta}$,
 associated to the braid $\gamma$, for some $0<\delta<1$. Recall
 that $f_{\gamma,\delta}(t,z)=(n.t, z(t)+h_\delta(z))$ where
 $h_\delta\colon \DD^2\to
\Int(\DD^2)$ is the homothety of ration $\delta$.

Consider $\varphi\in \Diff_0^1(\DD^2,\Int(\DD^2))$. By using
Corollary~\ref{c.parthyp}, one just needs to prove
Proposition~\ref{p.realization}  for a conjugate of $\varphi$ by
an orientation preserving diffeomorphism of $\DD^2$.

This allows us to assume that $\varphi(\DD^2)$ is contained in the
disk $\DD^2(\delta)$ of radius $\delta$ and that there is a
differentiable isotopy from $\varphi$ to the homothety $h_\delta$,
whose image remains contained in $\DD^2(\delta)$. More precisely,
there is a $C^1$-map $\Phi\colon \DD^2\times [-1,1]\to\Int(\DD^2)$
of the form $\Phi(x,t)=\varphi_t(x)$ where:
\begin{itemize}
\item for every $t\in [-1,1]$ one has $\varphi_t\in
\Diff_0^1(\DD^2,\Int(\DD^2))$, \item for every $t\in [-1,1]$ one
has $\varphi_t(\DD^2)\subset \DD^2(\delta)$, \item
$\varphi_0=\varphi$, and \item $\varphi_t=h_\delta$ if $|t|\geq
\frac 12$.

\end{itemize}
We denote $C=\max_{t\in[-1,1],z\in\DD^2} \|D_z(\varphi_t)\|$.

Let $\psi\colon [-\frac1{|n|},\frac1{|n|}]\to [-1,1]$ be a
diffeomorphism such that there is $0<\varepsilon <\frac 1{2|n|}$
with the following properties:
\begin{itemize}
\item $\psi(\frac 1n)=1$ and $\psi(-\frac1n)=-1$; \item for every
$t\in [-\frac1{|n|},\frac1{|n|}]$ one has $\left| \frac
d{dt}\psi(t)\right|>1$; \item $\left| \frac d{dt}\psi(t)\right|>
2C$ for every $t\in[-\varepsilon,\varepsilon]$; \item $\left|
\frac d{dt}\psi(t)\right|=|n|$ for $|t|\geq
\frac1{|n|}-\varepsilon$.
\end{itemize}

We define $f\colon S^1\times \DD^2\to \Int(S^1\times \DD^2)$ as
follows:
\begin{itemize}
\item $f(t,z)=(\psi(t), z(\frac{\psi(t)}{|n|})+ \varphi_{\frac
t\varepsilon}(z))$ if $t\in[-\varepsilon,\varepsilon]$,

\item $f(t,z)=(\psi(t),z(\frac{\psi(t)}{|n|})+ h_\delta(z))$ if
$|t|\in[\varepsilon,\frac1{|n|}]$,

\item $f(t,z)=(n.t,z(t)+h_\delta(z))$ if $t\notin[-\frac
1{|n|},\frac 1{|n|}]$.
\end{itemize}

Notice that $f(\{0\}\times \DD^2)\subset \{0\}\times \DD^2$, the
disk is normally hyperbolic and the restriction of $f$ to that
disk induces $\varphi$.  One concludes the proof of
Proposition~\ref{p.realization} by proving:

\begin{clai} The map $f$ defined above belongs to
$\cU_\gamma^{part.hyp}$.
\end{clai}
\begin{demo} We first notice that the image  $f(t,0)$ belongs to
the curve $\gamma(S^1)$; in other words $f(t,0)=\gamma(\tau_t)$,
where $t\mapsto\tau_t$ is a diffeomorphism of the circle. So the
image of $\{t\}\times \DD^2$ is contained in a disc of radius
$\delta$ in $\{n.\tau_t\}\times \DD^2$ centered at
$\gamma(\tau_t)$. As a consequence, if $r\neq s$ then
$f(\{r\}\times \DD^2)\cap f(\{s\}\times \DD^2)=\emptyset$. One
deduces that $f$ is injective, hence is a diffeomorphism from
$S^1\times \DD^2$ onto its image contained in
$f_{\gamma,\delta}(S^1\times \DD^2)$. One deduces that $f$ belongs
to $\cU_\gamma$.

In order to get the partial hyperbolicity, we will verify that $f$
satisfies the hypotheses of Lemma~\ref{l.parthyp} in each of the
possible expressions. We first notice that $f$ keeps invariant the
trivial foliation of $S^1\times \DD^2$ by the disks $\{t\}\times
\DD^2$. It remains to get the control of the derivative of $f$.

The map $f$ coincides with $f_{\gamma,\delta}$ out of $[-\frac
1{|n|},\frac1{|n|}]\times \DD^2$, giving the condition in this
region. On $\left([-\frac1{|n|},-\frac 1{|n|}+\varepsilon]\cup [
\frac1{|n|}-\varepsilon,\frac1{|n|}]\right)\times \DD^2$, one
notices that the derivative of the restriction of $f$ to each disk
$\{t\}\times \DD^2$ is the homothety of ratio $0<\delta<1$; hence
the conclusion holds because  $\left|\frac d{dt}
\psi(t)\right|>1$. Finally, for $t\in[\-\varepsilon,\varepsilon]$
the derivative of the restriction of $f$ to the disk $\{t\}\times
\DD^2$ is bounded by the constant $C$, and $\left|\frac d{dt}
\psi(t)\right|>C$, by assumption.
\end{demo}

\end{demo}

\section{An attracting solid torus $S^1\times
\DD^2$ without attractors.}\label{ss.torus}

Our main results are consequences of a construction in the solid
torus $S^1\times \DD^2$, that we explain in this section.

\subsection{Plykin attractors on normally hyperbolic disks, for solenoid maps}

Recall that $\cP_0$ is the open set of structurally stable
diffeomorphisms in $\Diff^1_0(\DD^2,\Int(\DD^2))$, defined at
Section~\ref{ss.plykyn}, whose non-wandering set consists exactly
in the union of a non-trivial hyperbolic attractor (a Plykin
attractor) and a finite set of periodic sources.

Given any braid $\gamma$ with $|n_\gamma|\geq 2$, let
$\cU_{\gamma}^{Ply}$ denote  the set of diffeomorphisms
$f\in\cU_\gamma^{part.hyp}$ such that $f$ leaves positively
invariant a normally hyperbolic essential disk $D_f$, and  such
that the restriction $\phi_f$ of $f$ to $D_f$ is $C^1$-conjugate
to an element $\phi\in\cP_0$.

As a  corollary of Proposition~\ref{p.realization} one gets:

\begin{coro}\label{c.realisation} Given any braid $\gamma$ with $|n_\gamma|\geq 2$, the set $\cU_{\gamma}^{Ply}$ is a non-empty $C^1$-open subset of
$\cU_\gamma^{part.hyp}$.
\end{coro}
\begin{demo}Proposition~\ref{p.realization} implies that
$\cU_{\gamma}^{Ply}$ is non-empty. It is open because the disk
$D_f$ is normally hyperbolic, hence persists by perturbation and
vary $C^1$-continuously with $f$; hence the restriction $\phi_f$
varies $C^1$-continuously with $f$; one concludes by recalling
that $\cP_0$ is an open subset of $\Diff^1_0(\DD^2,\Int(\DD^2))$.
\end{demo}

Let  $\gamma\subset S^1\times \DD^2$ be a braid  with
$|n_\gamma|\geq 2$. Consider  $f\in\cU_{\gamma}^{Ply}$. By
definition of $\cU_{\gamma}^{Ply}$ and $\cP_0$, one has the
following properties:

\begin{itemize}
\item there are $\alpha_f>0$ and $\ell>0$ such that the cone field
$\cC_{\alpha_f}$ is strictly invariant by $Df^\ell$;

\item   the disk $D_f$ is  positively invariant and normally
hyperbolic; hence  the disk $D_f$ is transverse to the cone field
$\cC_{\alpha_f}$; \item  the restriction $\phi_f$ of $f$ to $D_f$
belongs to $\cP_0$; hence, the disk $D_f$ contains a Plykin
attractor $A_f$ of $\phi_f$; as the disk $D_f$ is normally
hyperbolic,  $A_f$ is a hyperbolic basic set for $f$;

\item the Plykin attractor $A_f$ contains a hyperbolic fixed point
$p_f= x_{\phi_f}$ such that $Det(D\phi_f(p_f))>1$; as a
consequence, the product of any two eigenvalues of $Df(p_f)$ has
modulus larger than $1$; hence, the point $p_f$ is a sectionally
expanding fixed point of $f$;

\item we denote by $\La_f$ the chain recurrence class of $A_f$; in
an equivalent way, $\La_f$  is the chain recurrence class of the
fixed point $p_f$.
\end{itemize}

\subsection{Statement of our main result}

\begin{theo}\label{t.torus}
Given any braid $\gamma$ with $|n_\gamma|\geq 2$,
\begin{enumerate}
\item\label{i.unique} For every $f\in\cU_{\gamma}^{Ply}$, the
chain recurrence class $\La_f$ is the unique quasi attractor of
$f$. \item  We denote

$$\cU_{wild,\gamma}=\{f\in\cU_{\gamma}^{Ply},
\La_f\cap\overline{\{\mbox{sources of } f \}}\neq \emptyset\}.$$

In particular, $\La_f$ is not an attractor for $f\in
\cU_{wild,\gamma}$.

Then, for every $r\geq 1$,  the subset $\cU_{wild,\gamma}$ is
residual in $\cU_\gamma^{Ply}$ for the $C^r$ topology.
\end{enumerate}
\end{theo}

According to \cite{BC}, for $C^1$-generic diffeomorphisms, the
$\omega$-limit set $\omega(x)$ of any generic point $x$ of the
manifold is a quasi attractor. Hence the item~\ref{i.unique} of
Theorem~\ref{t.torus} implies:
\begin{coro} \label{c.torus} There is a $C^1$-residual subset of $\cU_{\gamma}^{Ply}$ of
diffeomorphisms $f$ for which the basin of $\La_f$ is residual in
$S^1\times \DD^2$.
\end{coro}

We don't know if Corollary~\ref{c.torus} holds for $C^r$-topology,
$r>1$. However, we think that it is possible to prove:
\begin{conj}
There is a $C^2$-open subset of $\cU_{\gamma}^{Ply}$ of
diffeomorphisms for which $\La_f$ carries an SRB-measure whose
basin has total Lebesgue measure in $S^1\times \DD^2$.
\end{conj}

\subsection{First step of the proof of Theorem~\ref{t.torus}: uniqueness of the quasi attractor}
\begin{demo}
Let $U\subset S^1\times \DD^2$ be an open attracting region of
$f$: $f(\overline U)\subset U$. Consider a segment $\sigma\subset
U$ which is tangent to $\cC_{\alpha_f}$. As $Df^\ell$ leaves
strictly invariant the cone field $\cC_{\alpha_f}$ and expands the
vectors in that cone field,  the forward iterates
$f^{n\ell}(\sigma)$, $n>0$, remain tangent to $\cC_{\alpha_f}$ and
their length tends to $\infty$. One deduces that there is $n>0$
such that $f^n(\sigma)\cap D_f\neq \emptyset$. Hence $f^n(U)\cap
D_f$ contains a non-empty open set. By definition of $\cP_0$ the
basin of the Plykin attractor $A_f$ of $\phi_f$ is a dense open
subset of $D_f$.  As a consequence, $f^n(U)$ contains a point $x$
in this basin. So $\omega(x,f)\subset A_f$. However
$\omega(x)\subset U$ because $U$ is by definition an attracting
region. So $U\cap A_f\neq\emptyset$. As $A_f$ is transitive and
$U$ is an attracting region, this implies $A_f\subset U$. In the
same way, the chain recurrence class $\La_f$ of $A_f$ is contained
in $U$.

Recall that a quasi attractor is a chain recurrence class which is
the intersection of a decreasing sequence of attracting regions.
This implies that every quasi attractor of $f$ contains $\La_f$,
hence is equal to $\La_f$.

On the other hand, as $S^1\times \DD^2$ is an attracting region,
it contains at least one quasi attractor. This concludes the
proof.
\end{demo}

\subsection{Robust homoclinic tangencies}

Now our construction consists in proving:
\begin{prop}\label{p.robusttangency}
For every $f\in \cU_\gamma^{Ply}$  and   every point $x$ of the
hyperbolic basic set $A_f$, there is $y\in A_f$ such that $W^u(x)$
and $W^s(y)$ meet tangentially  at one point.
\end{prop}
\begin{demo}[Idea of the proof of Proposition~\ref{p.robusttangency}]
Proposition~\ref{p.robusttangency} is completely analogous to
\cite[Proposition 3.1]{As}. We just recall the ideas of the proof
for completeness.

The hyperbolic set $A_f$ is a hyperbolic attractor of $\phi_f$.
Furthermore, by hypothesis, the basin $W^s(A_f)$ contains  the
whole disk $D_f$, punctured by the finite set $R_f$ of repelling
periodic points (contained in $f(D_f)$). Hence $D_f\setminus R_f$
is foliated by the stable manifolds of the point in $A_f$. Let us
denote $\cF^s$ this foliation.

On the other hand, as $D_f$ is an essential disk, transverse to
the cone field $\cC_{\alpha_f}$, for every circle $S^1\times\{z\}$
the image $f(S^1\times \{z\})$ cuts transversely $D_f$ in exactly
$|n_\gamma|>1$ points, (always with the same orientation). In
particular, $f(S^1\times \DD^2)$ cuts $D_f$ in exactly
$|n_\gamma|$ connected components, and one of them is $f(D_f)$.
Let $\De(f)$ be another component. Notice that $\cF^s$ induces a
foliation of the disk  $\De(f)$, already denoted by $\cF^s$.

Recall that $A_f$ is a lamination whose leaves are the unstable
leaves for $\phi_f$ of the points $z\in A_f$; these leaves are
tangent to the center-unstable direction of $A_f$ considered as a
basic set for $f$, and we denote them $L^c(z)$.

The unstable leaves of the points of $A_f$ are $C^1$ surfaces.
More precisely, for every point $z\in A_f$, the unstable manifold
$W^u(z)$ for $f$ is the union of all the strong unstable leaves
$L^{uu}(x)$  for $x\in L^c(z)$.

Each strong unstable leaf is a curve tangent to the cone field
$\cC_{\alpha_f}$, contained in $f(S^1\times \DD^2)$ and of
infinite length. In particular, every sufficiently large segment
of strong unstable leaf cuts the disk $\De(f)$. We endow the
strong unstable leaves with the orientation induced by the
orientation of the circle $S^1$. Hence, for any point $z\in A_f$
one has a well defined point $h(z)\in \De(f)$ which is the first
intersection point of $L^{uu}(z)$ with $\De(f)$. Notice that the
map $h$ is continuous.

Now  $\cL_f= h(A_f)$ is a regular $1$-dimensional compact
lamination contained in $\De(f)$.  Moreover, the  leaves are $C^1$
curves varying $C^1$-continuously, because they are obtained as
the (transverse) intersection of $\De(f)$ with the unstable
manifolds of the points $z\in A_f$.

Given any compact $1$-dimensional lamination by uniformly $C^1$
curves of a  $2$-disk endowed with a non-singular foliation, every
leaf of the lamination admits tangency points with the foliation.
So every leaf of the lamination $\cL_f$ admits tangency points
with $\cF^s$, ending the proof.
\end{demo}

\subsection{Proof of Theorem~\ref{t.torus}}

For proving Theorem~\ref{t.torus} we will show:
\begin{prop}\label{p.torus}
Given any braid $\gamma$ with $|n_\gamma|\geq 2$, and any
$\varepsilon>0$, the set

$$\cU_{n,\gamma}=\{f\in\cU_{\gamma}^{Ply}, \exists\ q_{n,f} \mbox{ hyperbolic periodic source, } d(p_f,q_{n,f})<\frac 1n\}$$
is open for the $C^1$ topology  and is dense in
$\cU_\gamma^{Ply}\cap {\rm Diff}^r(S^1\times\DD^2,\Int(S^1\times
\DD^2))$ for the $C^r$ topology, for every $r\geq 1$.
\end{prop}
Notice that $\bigcap_{n\in\NN^*} \cU_{n,\gamma}\subset
\cU_{wild,\gamma}$. Hence Proposition~\ref{p.torus} implies that
$\cU_{wild,\gamma}$ is residual in $\cU_\gamma^{Ply}$ for the
$C^r$ topology, for any $r\geq 1$, ending the proof of
Theorem~\ref{t.torus}.

The fact that $\cU_{n,\gamma}$ is $C^1$-open is a simple
consequence of the continuous dependence of hyperbolic periodic
points, for the $C^1$ topology. The difficulty is to prove the
$C^r$-density. As the set of $C^{r+1}$ diffeomorphisms is dense in
the set of $C^r$ diffeomorphisms for the $C^r$ topology,  the
$C^r$-density  of $\cU_{n,\gamma}$ is implied by the
$C^{r+1}$-density: hence it is enough to prove the $C^r$-density
of $\cU_{n,\gamma}$ for $r$ large enough.

However, the $C^1$-density can be proved by an argument of
different nature, involving specific $C^1$-perturbations lemmas.
We present both argument in the next sections.

\subsection{$C^1$-density of $\cU_{n,\gamma}$ }
Recall that $\cU_\gamma^{Ply}$ is contained in
$\cU_\gamma^{part.hyp}$. Hence every $f\in \cU_\gamma^{Ply}$ is
partially hyperbolic on $S^1\times\DD^2$: there is a dominated
splitting $T_x(S^1\times \DD^2)=E^{cs}(x)\oplus_{_<} E^u(x)$, for
$x\in \bigcap_{n\in\ZZ} f^n(S^1\times \DD^2)$, where $dim
E^{cs}=2$, $dim E^u=1$, and the vectors in $E^u$ are uniformly
expanded.

 For every $f\in \cU_\gamma^{Ply}$, consider the set $\Si_f$ of hyperbolic periodic saddle points which are homoclinically related with the point $p_f$ (i.e. whose stable and unstable manifolds cut transversally the unstable and stable manifold of $p_f$, respectively). Let $\Si_{f,0}\subset \Si_f$ be the set  of saddle points $p\in\Si_f$  which are sectionally expanding; in other words, $p\in\Si_f$ belongs to $\Si_{f,0}$ if,
 $$ \left|Det\left(Df^{\pi(p)}|_{E^{cs}(p)}\right)\right|>1,$$
 where $\pi(p)$ is the period of $p$ and $Df^{\pi(p)}|_{E^{cs}(p)}$ is the restriction of the derivative at the period to the center stable bundle at $p$.

Recall that the point $p_f$ is sectionally expanding (i.e.
$p_f\in\Si_{f,0}$). A classical argument, formalized by using the
notion of transitions in \cite{BDP} and used by many authors,
implies that for every $f$ the set $\Si_{f,0}$ is dense in the
homoclinic class $H(p_f,f)$ of $p_f$ (i.e. the closure of
$\Sigma_f$). More precisely, there is a sequence
$p_{f,i}\in\Si_{f,0} $, $i\in\NN$, such that, for every $\delta
>0$ and for every $i$ large enough, the orbit of $p_{f,i}$ is
$\delta$-dense in $H(p_f,f)$. We denote by $\pi_i$ the period of
$p_{f,i}$.

  Now, according to \cite{BC}, for $C^1$-generic $f\in\cU_\gamma^{Ply}$, the homoclinic class $H(p_f,f)$ coincides with the chain recurrence class $\La_f$. As  a consequence one gets:
\begin{lemm}For $C^1$-generic $f\in\cU_\gamma^{Ply}$, the closure of $\Si_{f,0}$ contains $\La_f$.
\end{lemm}
 According to Proposition~\ref{p.robusttangency}, for any $f\in\cU_\gamma^{Ply}$ the chain recurrence class $\La_f$ contains a tangency point $q_f$ of $W^u(p_f)$ with $W^s(A_f)$. One deduces:

\begin{lemm} For every $f\in\cU_\gamma^{Ply}$ the $2$-dimensional bundle $E^{cs}$ does not admit any dominated splitting along $\La_f$.
\end{lemm}
\begin{demo}We argue by contradiction, assuming that there is a dominated splitting $E^{cs}=E_1\oplus_{_<} E_2$ on $\La_f$: this splitting defines a dominated splitting $T_{\La_f}(S^1\times \DD^2)=E_1\oplus_{<}E_2\oplus_{_<} E^u$ on $\La_f$.  Then the stable manifold  $W^u(p_f)$  is tangent to $E_2\oplus E^u$ . Furthermore, for every $x\in A_f$ and every $y\in W^s(x)\cap \La_f$, the stable manifold $W^s(x)$ is tangent to $E_1(y)$ at $y$. This prevents $W^u(p_f)$ to have a tangency point with $W^s(x)$ for $x\in A_f$, hence contradicts Proposition~\ref{p.robusttangency}.
\end{demo}

As a direct corollary one gets:

\begin{coro}For every $C^1$-generic $f\in\cU_\gamma^{Ply}$, the $2$-dimensional bundle $E^{cs}$ does not admit any dominated splitting along $\Si_{f,0}$.
\end{coro}

Now, an argument of Ma\~n\'e in \cite{Ma}  (see also \cite{BDP})
shows that, for every $\varepsilon>0$ and every $i$ large enough,
there is an $\varepsilon$-$C^1$-perturbation $g_i\in
\cU_{\gamma}^{Ply}$ of $f$ which coincides with $f$ on the orbit
of $p_{f,i}$ and out an arbitrarily small neighborhood of this
orbit, and such that the (real or complex) eigenvalues of
$Dg_i^{\pi_i}(p_{f,i})$ corresponding to the center-stable bundle
$E^{cs}$ have the same modulus; furthermore, as $p_{f,i}$ was
sectionally expanding for $f$, this modulus can be taken larger
than $1$; as the eigenvalue corresponding to the unstable bundle
is also larger than $1$ one gets that the orbit of  $p_{f,i}$ is a
hyperbolic source for $g_i$. Hence, choosing $\varepsilon>0$ small
enough (so that the continuation $p_{g_i}$ of $p_f$ remains
arbitrarily close to $p_f$) and $i$ large enough (so that the
orbit of $p_{f,i}$ is passing arbitrarily close to $p_f$) one gets
$g_i\in\cU_{n,\gamma}$, ending the proof of the density of
$\cU_{n,\gamma}$ in $\cU_\gamma^{Ply}$ for the $C^1$ topology.

\subsection{$C^r$-density of $\cU_{n,\gamma}$ for $r\geq 2$}
We consider now $\cU_\gamma^{r,Ply}=\cU_\gamma^{Ply}\cap
\Diff^r(S^1\times \DD^2, \Int(S^1\times \DD^2))$ endowed with the
$C^r$-topology, for $r\geq 2$.

According to  Proposition~\ref{p.robusttangency} for every $f\in
\cU_\gamma^{r,Ply}$, the unstable manifold $W^u(p_f)$ presents a
tangency point $q_f$ with the stable manifold of a point $z_f\in
A_f$. Notice that $W^u(p_f)$ and $W^s(z_f)$ are $C^r$-immersed
submanifold, and $r\geq 2$. By performing an arbitrarily small
$C^r$ perturbation of $f$, one may assume that the tangency point
$q_f$ is a quadratic tangency point.

Then, for every $g$ in a small $C^2$-neighborhood $\cV$ of $f$ the
tangency point $q_f$ of $W^u(p_f)$ with the stable foliation of
$A_f$ has a unique continuation $q_g$, quadratic tangency point of
$W^s(p_g)$ with the stable foliation of $A_g$. This tangency point
varies continuously with $g$.

Notice that the positive orbit of $q_f$ is contained in the
invariant normally hyperbolic disk $D_f$ containing $A_f$. The
negative orbit of $q_f$ is not contained in $D_f$: by
construction, $q_f$ belongs to the lamination $\cL_f=h(A_f)$,
hence, is  the first return map on $D_f$ of the strong unstable
leaf of a point $y_f\in A_f$ (i.e. $q_f=h(y_f)$); so for $n>0$
large $f^{-n}(q_f)$ is a point contained in the local strong
unstable leaf of $f^{-n}(y_f)\in A_f\subset D_f$. So, one can
perform small $C^r$-perturbation of $f$ in a neighborhood of
$f^{-n}(q_f)$ without modifying the restriction of $f$ to the disk
$D_f$, hence without modifying the stable foliation $\cF_f$ of
$A_f$ in $D_f$. So we get:
\begin{lemm} There is a $C^r$ arc $\{f_t\}, t\in[0,1]$ of $C^r$ diffeomorphisms $f_t\in\cV$ such that:
\begin{itemize}
\item $f_0=f$; \item  for every $t\in[0,1]$, $f_t$ coincides with
$f$ on the disk $D_f$ (in particular  the stable foliation
$\cF_{f_t}$ of $A_{f_t}$ is $\cF_f$); \item the tangency point
$q_t=q_{f_t}$ defines an arc transverse to the stable foliation
$\cF_f$.
\end{itemize}
\end{lemm}

Recall that $A_f$ is a (transitive) hyperbolic attractor for the
restriction on $f$ to $D_f$, and the fixed point $p_f$ belongs to
$A_f$. Hence the stable manifold of $p_f$ is a dense leaf of the
foliation $A_f$. As a consequence one gets:
\begin{coro}There is a sequence $t_n>0$ tending to $0$ such that, for every $n\in \NN$, the point $q_{t_n}$ is a quadratic tangency point of the stable manifold of $p_f$ with the unstable manifold of $p_f$, and $p_f$ is a hyperbolic sectionally expanding point of $f_{t_n}$.
\end{coro}

Hence $g=f_{t_n}$ is an arbitrarily small $C^r$-perturbation of
$f$ having a quadratic homoclinic tangency point associated to a
sectionally expanding fixed point $p_g=p_f$. This situation has
been studied in \cite{PV}:

\begin{theo} \cite{PV}  If $\{g_s\}_{s\in [0,1]}$ is a generic arc of $C^r$ diffeomorphisms ($r\ge 2$), and  there is a periodic hyperbolic point $p$ of $g_0$ which is sectionally expanding, and such that $W^s(p,g_0)\cap W^u(p,{g_0})$ contains a quadratic tangency point $q$. Then there are a sequence $s_i$ converging to $0$ and periodic sources $q_i$ of $g_{s_i}$ converging to  $q$.
\end{theo}
Notice that, for large $i$,  the orbits of the periodic sources
$q_i$ are passing arbitrarily close to the point $p$. As a
consequence, for $i$ large the diffeomorphism $g_{s_i}$ belongs to
$\cU_{n,\gamma}$, and is an arbitrarily $C^r$-small perturbation
of $g$ which is an arbitrarily $C^r$-small perturbation of $f$.
This proves the $C^r$-density of $\cU_{n,\gamma}$ in
$\cU_\gamma^{Ply}$, ending the proof of Proposition~\ref{p.torus}.

\section{Non-existence of attractors for diffeomorphisms}\label{ss.dim3}

\subsection{An attracting ball $B^3$ without attractors}
Theorem~\ref{Thm:3Dnoattractor} is obtained from
Theorem~\ref{t.torus} by building locally generic diffeomorphisms
of an attracting ball $B^3$ without topological attractors and
with a unique quasi attractor:

\begin{theo}\label{t.ball}
There is a non-empty $C^1$-open subset $\cU\subset {\rm Diff}^1(
\DD^3, \Int(\DD^3))$ and, for $f\in \cU$, a hyperbolic periodic
point $p_f$ varying continuously with $f$ such that:
\begin{enumerate}
\item the diffeomorphism $f$ is $C^1$ conjugated  with the
homothety $z\mapsto\frac12z$ in a neighborhood of the sphere
$\partial \DD^3$; \item for every $f\in\cU$, the chain recurrence
class $\La_f=C(p_f)$ is the unique quasi attractor of $f$; \item
for every $r\geq 1$, the subset
$$\cU_{wild}=\{f\in\cU, \La_f\cap\overline{\{\mbox{sources of } f \}}\neq \emptyset\}$$
is residual for the $C^r$ topology.
\end{enumerate}
\end{theo}

Next lemma can be easily proved by using the same kind of
perturbations used for the \emph{derived from Anosov}
diffeomorphisms in \cite{Sm}. We leaves the details of the
construction to the reader.
\begin{lemm}\label{l.DA}Let $f\colon S^1\times \DD^2\to \Int(S^1\times \DD^2)$ be a solenoid map such that $\bigcap_{n\in\NN} f^n(S^1 \times \DD^2)$ is a hyperbolic attractor.  Then, there is $g$ isotopic to $f$, which coincides with $f$ in a neighborhood of the boundary $\partial (S^1\times \DD^2)$, and such that the chain recurrent set in $S^1\times \DD^2$ consists in exactly one fixed hyperbolic sink $\omega$ and a hyperbolic basic set of saddle type (i.e. neither attracting nor repelling). Moreover, if $f$ is orientation preserving,  one may require that  the derivative $Dg(\omega)$ is the homothety of ratio $\frac 12$.
\end{lemm}

\begin{demo}[Proof of Theorem~\ref{t.ball}]
According to \cite{Gi} there is a diffeomorphism $f_0$ of the $3$
sphere $S^3$ admitting a  torus $T$ with the following properties:
\begin{itemize}
\item the torus $T$ bounds two solid tori $\De_1$ and $\De_2$;
\item $f_0(\De_1)$ is contained in the interior of $\De_1$ and the
restriction $f_0|_{\De_1}$ is a hyperbolic Smale-solenoid
attractor corresponding to a $2$-braid $\gamma$; \item
$f_0^{-1}(\De_2)$ is contained in the interior of $\De_2$ and the
restriction $f_0^{-1}|_{\De_2}$ is a hyperbolic Smale-solenoid
attractor corresponding to a $2$-braid $\gamma$.
\end{itemize}

We now modify $f_0$  by surgery in both solid tori $\De_1$ and
$\De_2$, in order to get a diffeomorphism $f_1$ with the following
properties:

\begin{itemize}
\item $f_1$ coincides with $f_0$ in the neighborhood of the torus
$T$; as a consequence  $f_1(\De_1)\subset \Int(\De_1)$ and
$f_1^{-1}(\De_2)\subset \Int(\De_2)$; \item the restriction of
$f_1$ to the solid torus $\De_1$ belongs to the $C^1$-open set
$f\in \cU_\gamma^{Ply}$; \item the intersection of the chain
recurrent set $\cR(f_1)$ with $\De_2$ consists exactly in a
hyperbolic fix source $\alpha_1$ and a non-trivial hyperbolic set
$K_1$ of saddle type (this is obtained by applying
Lemma~\ref{l.DA} to the restriction of $f^{-1}$ to the solid torus
$f(\De_2)$).
\end{itemize}

Now one removes from $S^3$ the interior of a small ball $B$
centered at $\alpha_1$. Then $\cB=S^3\setminus \Int (B)$ is a
compact ball diffeomorphic to $\DD^3$. Furthermore $f_1(\cB)$ is
contained in the interior of $\cB$. Now there is a $C^1$
neighborhood $\cU$ of $f_1$ such that every $f\in\cU$ satisfies
the following properties:
\begin{itemize}
\item there is a diffeomorphism $\varphi\colon \cB\to \DD^3$ such
that $\varphi f\varphi^{-1}\colon \DD^3\to \DD^3$ coincides with
the homothety $z\mapsto\frac 12 z$ in a neighborhood of
$S^2=\partial \DD^3$; \item the image of the solid torus
$\De_1\subset \cB$
 is contained in its interior and the restriction $f|_{\De_1}$ belongs to $\cU_\gamma^{Ply}$; one denotes by $\La_f$ the unique quasi attractor of $f$ contained in $\De_1$, and by $p_f$ the hyperbolic sectionally expanding saddle point in $\La_f$ associated to $f|_{\De_1}\in\cU_\gamma^{Ply}$;
\item the intersection of the  chain recurrent set $\cR(f)$ with
$\cB\setminus \Int(\De_1)$ is a hyperbolic basic set of saddle
type.
\end{itemize}

 One concludes by noticing that $C^r$-generic diffeomorphisms $f\in\cU$ induce by restriction on $\De_1$ $C^r$-generic diffeomorphisms in $\cU_\gamma^{Ply}$; as a consequence, there is a sequence of hyperbolic sources converging to a point in $\La_f$, preventing $\La_f$ to be an attractor.
\end{demo}

\subsection{End of the proof of Theorem~\ref{Thm:3Dnoattractor}}

For getting Theorem~\ref{Thm:3Dnoattractor} one considers the time
one map of the flow of a gradient vector field of a Morse function
on $M$. Then one replaces the diffeomorphism in a neighborhood of
each sink by a diffeomorphism in the open set  $\cU$ built in
Theorem~\ref{t.ball}.

\begin{rema} Let $M$ be a compact orientable $3$-manifold.  Using
the fact that $M$ admits a Heegaard splitting in two handelbodies,
one easily verifies that $M$ admits a gradient like diffeomorphism
having a unique sink.  As a consequence, we can assume that $k=1$
in the statement of Theorem~\ref{Thm:3Dnoattractor}.
\end{rema}

\subsection{Non existence of attractors and repellers in higher
dimensions: proof of
Theorem~\ref{Thm:4Dnoattractorrepellers}}\label{ss.dim4}

Multiplying our construction in $B^3$  by a transverse contraction
allows us to get
\begin{lemm}\label{l.d-ball} Given any $d>3$, there is a non-empty $C^1$-open subset $\cU_d\subset {\rm Diff}^1( \DD^d,
\Int(\DD^d))$ such that  every $f\in \cU_d$ satisfies the
following properties:
\begin{enumerate}
\item the diffeomorphism $f$ is $C^1$ conjugated  with the
homothety $z\mapsto\frac12z$ in a neighborhood of the sphere
$\partial \DD^d$; \item the chain recurrent set of $f$ is
contained in a normally hyperbolic $3$-disc $D_f$; \item the
restriction $f|_{D_f}$ belongs to the open subset $\cU$ given by
Theorem~\ref{t.ball}; in particular
 for every $f\in\cU_d$, the chain recurrence class $\La_f$ of the fixed point $p_f$ is the unique quasi attractor of $f$.
\end{enumerate}
As a consequence, for every $r\geq 1$, the subset
$$\cU_{wild}=\{f\in\cU_d, \La_f\cap\overline{\{\mbox{sources of the restriction } f|_{D_f} \}}\neq \emptyset\}$$
is residual for the $C^r$ topology. Then $f\in\cU_{wild}$ has
neither attractors  nor repellers in $\DD^d$.

\end{lemm}

Given any manifold $M$ with $dim(M)>3$, one considers a
diffeomorphism $f_0$ which is the time one map of a Morse
function. Now, one builds a diffeomorphism $f_1$ obtained from
$f_0$ as follows:
\begin{itemize}
\item one replaces $f_0$, in a small ball centered to each sink,
by a diffeomorphism in the open set $\cU_d$ built at
Lemma~\ref{l.d-ball}; \item one replaces $f_0^{-1}$, in a small
ball centered to each source, by the inverse of  a diffeomorphism
in the open set $\cU_d$ built at Lemma~\ref{l.d-ball}.
\end{itemize}

Now the open set announced in
Theorem~\ref{Thm:4Dnoattractorrepellers} is obtained by
considering a small neighborhood of the diffeomorphism $f_1$
above.

\section{Singular flows: proof of Theorem~\ref{Thm:4Dvectorfield}\label{ss.vector}}

Our example for flow is very similar to the examples built for
Theorem~\ref{t.torus}, so that we will just sketch the
construction.

We consider an open set $\cU$ of vector fields on  $\RR^4$, such
that every $X\in\cU$ satisfies the following properties:
\begin{itemize}

\item the vector field $X$ admits  a transverse cross section
$\Si$ diffeomorphic to a solid torus $S^1\times \DD^2$;

\item the vector field $X$ has a unique singular point $0_X$ which
is a saddle with $dim (W^s(0_X))=3$;  the eigenvalues of the
derivative $D_{0_X}X$ are
$$\lambda_1<\lambda_2<\lambda_3<0<\lambda_4,$$
with $\lambda_4+\lambda_1>0;$

\item there is an essential disc $D_0\subset \Si$, transverse to
all the circles $S^1\times\{z\}, z\in \DD^2$,  and contained in
the local stable manifold of the saddle point $0_X$;

\item the first return map on $\Si$ is well defined on
$\Si\setminus D_0$ and the image is contained in the interior of
$\Si$; we denote it $P\colon \Si\setminus D_0\to \Int(\Si)$;

\item the first return map $P$ leaves invariant a splitting
$T\Si=E^{cs}\oplus E^u$ which is a dominated splitting with $\dim
E^{cs}=2$ and $\dim E^u=1$; moreover, $E^u$ is transverse to the
discs $\{t\}\times \DD^2$, $t\in S^1$;

\item the bundle $E^u$ is uniformly expanding by a factor larger
than $3$; more precisely, given any non-zero vector $u$ tangent to
$E^u(x)$, $x\in\Si$, denote $u=u_h+u_v$ where $u_h$ is tangent to
the $S^1$ fiber through $x$ and $u_v$ is tangent to the $\DD^2$
fiber through $x$;

    assume $x\in \Si\setminus D_0$ and let $w=D_xP(u)=w_h+w_v$; then we require:
    $$|w_h|>3|u_h|;$$

\item there is an essential disc $D_1\subset \Si\setminus D_0$,
invariant by $P$ (i.e. $P(D_1)\subset \Int(D_1)$), normally
hyperbolic, and such that the restriction $P|_{D_1}$ is smoothly
conjugated to an element of the open set $\cP_0$  of structurally
stable diffeomorphisms in $\Diff^1_0(\DD^2,\Int(\DD^2))$, defined
at Section~\ref{ss.plykyn}; in particular, the chain recurrent set
of   $P|_{D_1}$ consists in a Plykin attractor $A_X$ and finitely
many repelling points, and the Plykin attractor  $A_X$ contains a
fixed point $p_X$ which is sectionally expanding;

\item there is an essential disc $D_2\subset \Si\setminus D_0$,
invariant by $P$, normally hyperbolic, and such that the
restriction $P|_{D_2}$ has a unique fixed point $q_X$; the disc
$D_2$ is contained in the stable manifold of $q_X$ (for the map
$P$); finally, the derivative $D_{q_X}(P)$ of $P$ at $q_X$ has a
complex (non-real) eigenvalue, corresponding to the tangent space
$T_{q_X}D_2$.
\end{itemize}

It is not hard to build a non-empty open set $\cU$ of vector
fields satisfying all the properties above (see also \cite{BLY}
which contains the details of  a analogous construction).

As in the proof of Theorem~\ref{t.torus}, one verifies that, for
any open subset $O\subset \Si$ there is $n>0$ such that $f^n(O)$
meets $D_0$, $D_1$ and $D_2$: this implies that every attracting
region for $X$ which meets $\Si$ contains the singular point
$0_X$, the Plykin attractor $A_X$ (and hence its orbits by the
flow of $X$) and the orbit $\gamma_X$ of the point  $q_X$. Hence
there is a unique quasi attractor $\La_X$ for the orbits of $X$
through $\Si$  and this quasi attractor contains $0_X$ , $A_X$ and
$\gamma_X$.  An analogous argument shows that, for every
$X\in\cU$,  the invariant manifolds of $A_X$ for $P$ present a
tangency point. This implies that $C^r$-generic pathes in $\cU$
unfold generic homoclinic bifurcations associated to $p_X$,
implying that, for $C^r$-generic $X\in \cU$ the quasi attractor
$\La_X$ is accumulated by periodic sources, which prevents $\La_X$
to be an attractor.

One concludes the proof of Theorem~\ref{Thm:4Dvectorfield} by
proving
\begin{lemm} For any $X\in\cU$, the tangent flow of $X$ on $\La_X$ does not admit any dominated splitting.
\end{lemm}
\begin{demo} Assume that there is a dominated splitting $TM|_{\La_X}=E\oplus_{_<} F$, for the tangent flow of $X$. This dominated splitting induces on $\Si\cap \La_X$ a dominated splitting $T\Si|_{\Si\cap \La_X}=E_\Si\oplus F_\Si$ invariant by $P$ (just consider $E_\Si= (E+\RR X)\cap T\Si$ and $F_\Si=(F+\RR X)\cap T\Si$).

The fact that $q_X$ belongs to $\La_X\cap\Si$  implies that $\dim
E_\Si=2$. One deduces that $E_\Si=E^{cs}$ and $F_\Si=E^u$. As a
consequence one gets two possibilities for the splitting $T_xM =
E(x)\oplus F(x)$ at $x\in\La_X\cap\Si$:
\begin{itemize}
\item either $E=E^{cs}\oplus \RR X$ and $F\subsetneq E^u\oplus \RR
X$, \item or $E\subsetneq E^{cs}\oplus \RR X$ and $F= E^u\oplus
\RR X$.
\end{itemize}

In the first case, $X$ is tangent to $E$ along $\La_X$. However,
$\La_X$ contains the unstable manifold of $0_X$ (a quasi attractor
always contains its unstable manifold). This manifold consists in
$0_X$ and  $2$ orbits of $X$. Hence $W^u(0_X)$ is tangent to $E$.
This implies that $E(0_X)$ contains the eigenspace corresponding
to $\lambda_4$, which contradicts the fact that $E$ is dominated
by $F$.

In the second case, $X$ is tangent to $F$ along $\La_X$.  However,
for $x\in A_X$, the space $E^{cs}(x)$ contains vectors tangent to
the hyperbolic attractor $A_X\subset D_1$, hence contains vectors
which are exponentially expanded by the derivative $DP^n$, for
$n\to +\infty$. This implies that the space $E(x)$ contains
vectors $u\in E(x)$ and a sequence of times $t_n\to +\infty$ such
that $X_{t_n}(x)=P^n(x)\in\Si$ and
$$\lim_{n\to+\infty}|(X_{t_n})_*(u)|=+\infty,$$ where $(X_t)_*$
denotes the derivative of the time $t$ of the flow of $X$.

On the other hand, $X(x)\in F(x)$ but $|(X_{t_n})_*(X(x))|$
remains bounded, contradicting the fact that $F$ dominates $E$.

Hence both cases lead to contradiction, ending the proof.

\end{demo}

\section{Changing the definition of attractors}\label{s.conclusion}
With the better understanding of the complexity of generic
dynamics, people tried the definition of attractors in order to
ensure their existence.

\subsection{Palis approach from the point of view of ergodic theory}
From the ergodic viewpoint, an attractor $\Lambda$ of $f$ should
satisfy the following
\begin{itemize}
\item{``\emph{indecomposable property}'': there is an ergodic
invariant probability measure $\mu$ such that ${\rm
supp}(\mu)=\Lambda$;} \item ``\emph{attracting property}": its
basin $B(\Lambda)$ has positive Lebesgue measure, where
$$x\in B(\Lambda) \Longleftrightarrow \lim_{n\to\infty}\frac{1}{n}\sum_{i=0}^{n-1}\delta_{f^i(x)}=\mu$$
(here $\delta_z$ stands the Dirac measure at the point $z$).

\end{itemize}

\begin{Conjecture}[Palis \cite{P1,P2,P3}]
There is a dense set ${\cal D}\subset  {\rm Diff}^r(M)$ such that
for any $f\in\cal D$, $f$ has only finitely many (ergodic)
attractors, and the union of the basins of attractors forms a full
Lebesgue measure set in $M$.
\end{Conjecture}

Palis completed his conjecture by continuity properties of the
basins of the attractors with respect to the diffeomorphism.

\subsection{A topological approach}

In \cite{H}, Hurley proved that, for generic homeomorphisms of a
compact manifold, the $\omega$-limit set of every generic point is
a quasi attractor, and he stated the conjecture
\begin{Conjecture}[Hurley]
For $C^r$-generic diffeomorphisms the $\omega$-limit sets of
generic points in $M$ are quasi attractors.
\end{Conjecture}
This conjecture has been proved in \cite{MP,BC} for the
$C^1$-topology and remains open in more regular topologies.

However, the information given by Hurley's conjecture is very
weak: every $C^0$-generic homeomorphism $h$ has uncountably many
quasi attractors, and the closure of the basin of each quasi
attractor has empty interior\footnote{The proof of this last fact
(the closure of each basin has empty interior) was found by the
first author of this present paper, writing this conclusion;
Hurley kindly wrote us that he did not notice this fact.}. In the
setting  of  the $C^1$ topology, \cite{BD2} shows that there are
locally generic diffeomorphisms having an uncountable family of
quasi attractors which are at the same time quasi repellers; in
particular the basin of each of them is reduced to the
quasi-attractor itself (which is a Cantor set).

Let us define a new notion of attractor which will allow us to
propose a new conjecture.
\begin{defi}\begin{itemize}
\item  A \emph{residual attractor} of a diffeomorphism $f$ is a
chain recurrence class admitting a neighborhood $U$ which is an
attracting region and such that the $\omega$-limit set of the
generic points in $U$ is $\La$.

\item A \emph{locally residual attractor} of a diffeomorphism $f$
is a chain recurrence class admitting an open set $U$  such that
the $\omega$-limit set of the generic points in $U$ is $\La$.
Notice here $U$ may not be a neighborhood of $\La$.
\end{itemize}
\end{defi}
\begin{rema}
\begin{itemize}

\item For $C^1$-generic diffeomorphisms, one can  deduce from
\cite{BC} that  the residual attractors are exactly the quasi
attractors which are isolated in the set of quasi attractors: they
admit a neighborhood disjoint from any other quasi attractor.
\item For $C^1$-generic diffeomorphisms, our notion of
\emph{locally residual attractor} coincides with the notion of
\emph{generic attractor} introduce by Milnor in \cite{Mi}. More
precisely, Milnor first defines a \emph{minimal attractor} for the
ergodic point of view: the basin has positive Lebesgue measure and
every proper subset's basin only has zero Lebesgue measure; then,
on the topological generic setting he writes: ``There is an
analogous concept of \emph{generic-attractor}. The definition will
be left to the reader''. Hence, a \emph{generic attractor} is an
invariant set whose basin is a locally residual set, and such that
every proper subset's basin is meager. As Hurley's conjecture is
proved for $C^1$-generic diffeomorphisms, Milnor's generic
attractors of a $C^1$-generic diffeomorphism  are its locally
residual attractors.

\end{itemize}
\end{rema}

The locally generic examples built in
Theorem~\ref{Thm:3Dnoattractor} have finitely many residual
attractors  and the union of their basin is a residual subset of
the whole manifold $M$. This motivates the following problem:

\begin{prob}\label{p.ga}
\begin{enumerate}
\item Is it true that  $C^r$-generic diffeomorphisms have at least
one (locally) residual attractor? \item For any $C^r$-generic
diffeomorphism, is it true that the $\omega$-limit set of every
generic point is a (locally) residual attractor?
\end{enumerate}
\end{prob}
(A positive answer to these questions is known for locally
residual attractors  of $C^1$-generic diffeomorphisms: see the
next section, devoted to the $C^1$-topology).

We would like to understand better these residual attractors, in
particular to understand if their are associated to periodic
orbits. Recall that the \emph{homoclinic class} of a periodic
orbit is the closure of the transverse intersection of its
invariant manifolds. It is an invariant compact set canonically
associated to the periodic orbit. \cite{BC} shows that, for
$C^1$-generic diffeomorphisms, the chain recurrence class of a
periodic orbit is its homoclinic class; as a consequence, isolated
chain recurrence classes of $C^1$-generic diffeomorphisms   are
homoclinic classes (in particular, this holds for topological
attractors). As we noticed above, the residual attractors are the
quasi attractors which are isolated in the set of quasi
attractors. It seems natural to ask:

\begin{prob}\label{p.gahc} Let $\La$ be a residual attractor of a generic diffeomorphism. Is $\La$ the homoclinic class of a periodic orbit?
\end{prob}

\subsection{Remarks on the $C^1$ topology}
For $C^1$ generic non-critical (i.e. far from homoclinic
tangencies) diffeomorphisms, \cite{Y} gave a positive answer to
Problems~\ref{p.ga} and \ref{p.gahc}   proving that every quasi
attractor is a homoclinic class. Since for $C^1$ generic
diffeomorphism, we can have only countably many homoclinic
classes, together with the results in \cite{MP,BC}, there is at
least one locally residual attractor; furthermore, the (countable)
union of the basins of the locally residual attractors is a
residual subset of the manifold.

In a forthcoming work, we can get more precise results for the
$C^1$ topology.
\begin{itemize}
\item{On the contrary of Theorem \ref{Thm:3Dnoattractor}, we can
prove that for two dimensional manifold $M^2$, there is a $C^1$
dense open set ${\cal U}\subset \Diff^1(M^2)$, such that for any
$f\in\cal U$, $f$ has a hyperbolic attractor;}

\item{As a complement of Theorem \ref{Thm:3Dnoattractor}, for any
compact three dimensional manifold $M^3$ without boundary, we can
construct a $C^1$ open set ${\cal U}\subset \Diff^1(M^3)$, such
that  $C^1$-generic $f\in\cU$ have neither attractors nor
repellers;\footnote{For the examples built in Theorem
\ref{Thm:3Dnoattractor}, there are infinitely many repellers. }}

\item{Together with S. Gan, we give a positive answer to
Problems~\ref{p.ga} and~\ref{p.gahc} in the setting of partially
hyperbolic splitting with $1$-dimensional center bundle
In these setting, we prove that for $C^1$ generic diffeomorphism,
every quasi attractor is a residual attractor.}

\end{itemize}

\noindent Christian Bonatti,

\noindent {\small Institut de Math\'ematiques de Bourgogne}

\noindent {\small Universit\'e de Bourgogne, Dijon 21004, FRANCE}

\noindent {\footnotesize{E-mail : bonatti@u-bourgogne.fr}}

\vskip 5mm

\noindent Ming Li,

 \noindent{\small Department of Mathematics}

\noindent{\small Beijing Institute of Technology, Beijing 100081,
P. R. China}

\noindent {\footnotesize {E-mail : liming@math.pku.edu.cn}}

\vskip 5mm \noindent  Dawei Yang,

\noindent {\small School of Mathematic Sciences}

\noindent {\small Peking University, Beijing 100871, P. R. China}

\noindent{\footnotesize{E-mail : yangdw@math.pku.edu.cn}}

\end{document}